\input amstex

\define\Pee{{\Bbb P}}
\define\Zee{{\Bbb Z}}
\define\Cee{{\Bbb C}}

\define\W{\widetilde}
\define\w{\tilde}

\def\exact#1#2#3{0\rightarrow{#1}\rightarrow{#2}\rightarrow{#3}\rightarrow0}

\define\Sym{\operatorname{Sym}}

\define\Supp{\operatorname{Supp}}

\define\Quot{\operatorname{Quot}}

\define\proof{\demo{Proof}}
\define\endproof{\qed\enddemo}

\define\theorem#1{\proclaim{Theorem #1}}
\define\lemma#1{\proclaim{Lemma #1}}
\define\proposition#1{\proclaim{Proposition #1}}
\define\corollary#1{\proclaim{Corollary #1}}
\define\claim#1{\proclaim{Claim #1}}

\define\section#1{\specialhead #1 \endspecialhead}
\define\ssection#1{\medskip\noindent{\bf #1}}

\def\mapup#1{\uparrow
    \rlap{$\vcenter{\hbox{$\scriptstyle#1$}}$}}

\loadbold
\documentstyle{amsppt}
\leftheadtext{Wei-ping Li and Zhenbo Qin}
\rightheadtext{Blowup formulae}

\topmatter
\title
On blowup formulae for the $S$-duality conjecture of Vafa and Witten
\endtitle
\author {Wei-ping Li$^1$ and Zhenbo Qin$^2$}
\endauthor
\address 
Department of Mathematics, HKUST, Clear Water Bay, Kowloon, Hong Kong
\endaddress
\email mawpli\@uxmail.ust.hk  \endemail
\address Department of Mathematics, Oklahoma State University, 
Stillwater, OK 74078, USA
\endaddress
\email  zq\@math.okstate.edu \endemail
\thanks
${}^1$Partially supported by the grant HKUST631/95P
\endthanks
\thanks 
${}^2$Partially supported by NSF grant DMS-9622564 and 
an Alfred P. Sloan Research Fellowship
\endthanks

\endtopmatter

\TagsOnRight
\document
\section{1. Introduction}

In \cite{V-W}, Vafa and Witten formulated some mathematical predictions 
about the Euler characteristics of instanton moduli spaces
derived from the $S$-duality conjecture in physics
(details will be given in section 3). 
 From these mathematical predictions, a blowup formula was proposed
based upon the work of Yoshioka \cite{Yos}. 
Roughly speaking, the blowup formula says that there exists 
a universal relation between the Euler characteristics of 
instanton moduli spaces for a smooth four manifold
and the Euler characteristics of instanton moduli spaces for 
the blowup of the smooth four manifold. The universal relation is
independent of the four manifold and related to some modular forms.
This blowup formula of Vafa and Witten is different from 
the ``conventional" blowup formula (see \cite{F-S}) 
in gauge theory which usually means a relation between 
the Donaldson invariants of a smooth four manifold and 
the Donaldson invariants of the blowup of the smooth four manifold.

In this paper, using the virtual Hodge polynomials introduced in \cite{D-K}, 
we shall verify the blowup formula of Vafa and Witten for 
the gauge group $SU(2)$ and its dual group $SO(3)$ 
when the underlying four manifold is an algebraic surface. 
Let $\phi: \W X \to X$ be the blowing-up of an algebraic surface $X$ 
at a point $x_0 \in X$, and $E$ be the exceptional divisor. 
For simplicity, we always assume that $X$ is simply connected. 
Fix a divisor $c_1$ on $X$, $\w c_1 = \phi^*c_1 - aE$ with $a = 0$ or $1$,
and an ample divisor $H$ on $X$ with odd $(H \cdot c_1)$. 
For an integer $n$, let $\frak M_{H}(c_1, n)$ be the moduli space
of Mumford-Takemoto $H$-stable rank-$2$ bundles with 
Chern classes $c_1$ and $n$, $\frak M^G_{H}(c_1, n)$ 
be the moduli space of Gieseker $H$-semistable rank-$2$ torsion-free sheaves 
with Chern classes $c_1$ and $n$, and $\frak M^U_{H}(c_1, n)$ be 
the Uhlenbeck compactification of $\frak M_{H}(c_1, n)$ from 
gauge theory \cite{Uhl, Do1, LiJ}. 
For $r \gg 0$, the divisors $H_r = r \cdot \phi^*H - E$ on 
$\W X$ are ample; moreover, all the moduli spaces $\frak M_{H_r}(\w c_1, n)$ 
(resp. $\frak M^G_{H_r}(\w c_1, n)$, $\frak M^U_{H_r}(\w c_1, n)$) 
can be naturally identified. So we shall use the notation
$\frak M_{H_\infty}(\w c_1, n)$ (resp. $\frak M^G_{H_\infty}(\w c_1, n)$, 
$\frak M^U_{H_\infty}(\w c_1, n)$) to denote the moduli space
$\frak M_{H_r}(\w c_1, n)$ (resp. $\frak M^G_{H_r}(\w c_1, n)$, 
$\frak M^U_{H_r}(\w c_1, n)$) with $r \gg 0$. 

For a complex variety $Y$ (not necessarily smooth, projective, 
or irreducible), let $e(Y; x, y)$ be the virtual Hodge polynomial of 
$Y$ introduced in \cite{D-K}. For a complex scheme $Y$, 
we define $e(Y; x, y) = e(Y_{\text{red}}; x, y)$.
It is known (see \cite{Ful}) that when the complex variety $Y$ 
(carrying the Zariski topology) is projective,
$e(Y; 1,1)$ is the topological Euler characteristic $\chi(Y)$ of $Y$ 
regarded as a complex space (carrying the usual topology). 
Thus for a complex scheme $Y$, $e(Y; 1,1)$ (which is $e(Y_{\text{red}}; 1,1)$
by our definition) is equal to 
the topological Euler characteristic $\chi(Y_{\text{red}})$. 

It is well-known that both the Gieseker moduli spaces and 
the Uhlenbeck compactification spaces are complex projective schemes
(in fact, the Uhlenbeck compactification spaces are complex varieties
by J. Li's definition \cite{LiJ}). Our first main result in this paper 
is the following blowup formula for the Gieseker moduli spaces.

\theorem{A} Let $(H \cdot c_1)$ be odd and $\w c_1 = \phi^*c_1 - aE$ 
with $a = 0$ or $1$. Then 
$$\sum_{n} e(\frak M_{H_\infty}^G(\w c_1, n); x, y) q^{n- {\w c_1^2 \over 4}} 
= (q^{1 \over 12} \cdot \W {\W Z}_a) \cdot 
\sum_{n} e(\frak M_{H}^G(c_1, n); x, y) q^{n- {c_1^2  \over 4}}$$
where $\W {\W Z}_a = \W {\W Z}_a(x, y, q)$
is a universal function of $x, y, q, a$ with 
$$\W {\W Z}_a(1, 1, q) = {\sum_{n \in \Zee} q^{(n+{a \over 2})^2} \over 
[q^{1 \over 24} \prod_{n \ge 1} (1 - q^n)]^2}.$$
In particular, setting $x=y=1$ yields the blowup formula
$$\sum_{n} \chi(\frak M_{H_\infty}^G(\w c_1, n)_{\text{red}}) 
   q^{n- {\w c_1^2 \over 4}} 
= {\sum_{n \in \Zee} q^{(n+{a \over 2})^2} \over \prod_{n \ge 1} (1 - q^n)^2} 
\cdot \sum_{n} \chi(\frak M_{H}^G(c_1, n)_{\text{red}}) 
   q^{n- {c_1^2  \over 4}}.$$  
\endproclaim

The universal function $\W {\W Z}_a(1, 1, q)$ was conjectured by
Vafa and Witten \cite{V-W}. The exponents of $q$ are written in the 
forms of ${n- {\w c_1^2 \over 4}}$ and ${n- {c_1^2  \over 4}}$ 
which stand for the instanton numbers. 
We remark that over finite fields $\Bbb F_q$, Yoshioka \cite{Yos} proved 
a blowup formula for the number of $\Bbb F_q$-rational points in 
the Gieseker moduli spaces over any algebraic surface.
Using the Weil Conjecture and the facts that the Gieseker moduli spaces
over $\Pee^2$ are projective and smooth and that the number of 
$\Bbb F_q$-rational points in these moduli spaces is a polynomial in $q$, 
Yoshioka proved the blowup formula for the Gieseker moduli spaces 
when $X = \Pee^2$ and $x = y$ (i.e. a blowup formula for 
the Poincare polynomials of the Gieseker moduli spaces over $\Pee^2$). 
Since the Gieseker moduli spaces are singular for general algebraic surfaces,
the Weil Conjecture can not be applied to general algebraic surfaces. 
So Yoshioka's blowup formula over finite fields can not be carried to 
the complex field. However, Yoshioka's results and methods provided 
us with valuable guidances.

For the Uhlenbeck compactifications, we have the following blowup formula.

\theorem{B} Let $(H \cdot c_1)$ be odd and $\w c_1 = \phi^*c_1 - aE$ 
with $a = 0$ or $1$. Assume that $\frak M_{H}(c_1, n)$ 
(respectively, $\frak M_{H_\infty}(\w c_1, n)$) is dense 
in the Gieseker moduli space $\frak M_{H}^G(c_1, n)$ 
(respectively, in $\frak M_{H_\infty}^G(\w c_1, n)$) for all $n$. Then 
$$\sum_{n} e(\frak M^U_{H_\infty}(\w c_1, n); x, y) 
q^{n- {\w c_1^2 \over 4}} = (q^{1 \over 12} \cdot \W Z_a) \cdot 
\sum_{n} e(\frak M^U_{H}(c_1, n); x, y) q^{n- {c_1^2 \over 4}}$$
where $\W Z_a = \W Z_a(x, y, q)$ is a universal function of $x, y, q, a$ with 
$$\W Z_a(1, 1, q) = {\sum_{n \in \Zee} q^{(n+{a \over 2})^2} \over 
q^{1 \over 12}(1 - q)}.$$
In particular, setting $x=y=1$ yields the blowup formula
$$\sum_{n} \chi(\frak M^U_{H_\infty}(\w c_1, n)) q^{n- {\w c_1^2 \over 4}} =  
{\sum_{n \in \Zee} q^{(n+{a \over 2})^2} \over (1 - q)} \cdot 
\sum_{n} \chi(\frak M^U_{H}(c_1, n)) q^{n- {c_1^2 \over 4}}.$$
\endproclaim

We remark that the assumption that $\frak M_{H}(c_1, n)$ is dense 
in $\frak M_{H}^G(c_1, n)$ determines the structure of $\frak M^U_{H}(c_1, n)$.
Moreover, by the Lemma 2.3 in \cite{F-Q}, if $(c_1 \cdot H)$ is odd 
and the anti-canonical divisor $(-K_X)$ is effective, 
then $\frak M_{H}(c_1, n)$ is dense in $\frak M_{H}^G(c_1, n)$.
Thus if both $(-K_X)$ and $(-K_{\W X})$ are effective,
then $\frak M_{H}(c_1, n)$ (respectively, $\frak M_{H_\infty}(\w c_1, n)$) 
is dense in $\frak M_{H}^G(c_1, n)$ 
(respectively, in $\frak M_{H_\infty}^G(\w c_1, n)$) for all $n$.
In particular, if $(-K_X)$ is positive and the point $x_0$ is contained 
in a curve $C \in |-K_X|$, then $\frak M_{H}(c_1, n)$ 
(respectively, $\frak M_{H_\infty}(\w c_1, n)$) 
is dense in the Gieseker moduli space $\frak M_{H}^G(c_1, n)$ 
(respectively, in $\frak M_{H_\infty}^G(\w c_1, n)$) for all $n$.

The idea for the proof of our theorems is as follows. 
First of all, we reduce the blowup formulae
to some universal relations among the virtual Hodge polynomials
of the moduli spaces $\frak M_{H_\infty}(\w c_1, n), \frak M_{H}(c_1, n)$,
and $\frak M_{H}^G(c_1, n)$. Then using standard techniques such as
taking double duals and elementary modifications, 
we stratify $\frak M_{H_\infty}(\w c_1, n)$ and $\frak M_{H}^G(c_1, n)$
into finite disjoint unions of locally closed subsets $W$. 
Roughly speaking, these subsets $W$ admit morphisms $f_W$ to 
$\frak M_{H}(c_1, k)$ for some $k \le n$. Moreover, 
these morphisms $f_W: W \to \text{Im}(f_W)$ are 
Zariski-locally trivial bundles, 
and the fibers are simple objects such as certain Grothendieck Quot-schemes
and certain open subsets in the projective spaces. 
Using the properties of virtual Hodge polynomials,
we are able to establish the universal relations among 
the virtual Hodge polynomials of $\frak M_{H_\infty}(\w c_1, n), 
\frak M_{H}(c_1, n)$, and $\frak M_{H}^G(c_1, n)$.

Notice that since the Euler characteristic is not a cobordism invariant, 
the Euler characteristics of the Uhlenbeck compactification spaces 
on smooth four manifolds are not necessarily smooth invariants. 
These Euler characteristics depend on the Riemannian metrics on 
the smooth four manifolds. Nevertheless, we think that our blowup formula 
for the Euler characteristics of the Uhlenbeck compactification spaces 
may hold for general smooth four manifolds in the following sense.
Fix a generic Riemannian metric $g$ on a smooth four manifold $X$.
On the connected sum $\overline X = X \# \overline \Pee^2$, 
following Donaldson's construction (see \cite{Do2}),
we take the metric $\overline g$ on $\overline X$ which is close to $rg$ on $X$ for $r \gg 0$ and to the Fubini-Study metric on $\overline \Pee^2$.
Then our blowup formula may hold for the Euler characteristics of 
the Uhlenbeck compactification spaces on the Riemannian four manifolds 
$(X, g)$ and $(\overline X, \overline g)$.

Our paper is organized as follows. In section 2, 
we review virtual Hodge polynomials and their basic properties. 
In section 3, we recall the $S$-duality conjecture in the form formulated 
by Vafa and Witten, and slightly modify the conjectured blowup formula.
In sections 4 and 5, we prove the blowup formula for 
the Uhlenbeck compactifications and the Gieseker moduli spaces respectively.

\medskip\noindent
{\bf Acknowledgments:} The authors thank William Banks, Jan Cheah,
Jim Cogdell, Robert Friedman, Sheldon Katz, Jason Levy, Weiping Li, 
and Yongbin Ruan for valuable helps and stimulating discussions.
The authors are very grateful to the referees for useful comments 
and suggestions which have greatly improved the exposition of the paper.  

\section{2. Virtual Hodge polynomials}

Virtual Hodge polynomials can be viewed as a convenient tool for 
computing the Hodge numbers of smooth projective varieties 
by reducing to computing those of simpler varieties.
Using Deligne's work \cite{Del} on mixed Hodge structures,
Danilov and Khovanskii \cite{D-K} introduced
virtual Hodge polynomials for a complex algebraic variety $Y$
(not necessarily smooth, projective, or irreducible).
The mixed Hodge structures which are defined on the cohomology 
$H_c^k(Y, \Bbb Q)$ with compact support coincide with the classical one 
if the complex variety $Y$ is projective and smooth. 
For each pair of integers $(s, t)$, define the virtual Hodge number
$$e^{s, t}(Y) = \sum_k (-1)^k h^{s, t}(H_c^k(Y, \Bbb Q)).$$ 
Then the virtual Hodge polynomials of $Y$ is defined by 
$$e(Y; x, y) = \sum_{s, t} e^{s, t}(Y) x^sy^t.  $$
Virtual Hodge polynomials satisfy the following properties
(\cite{D-K, Ful, Che}):

\roster
\item"{(2.1)}" When $Y$ is projective, $e(Y; 1, 1)$ is the same as 
the Euler characteristic $\chi(Y)$.
When $Y$ is projective and smooth, $e(Y; x, y)$ is the usual Hodge polynomial. 
\item"{(2.2)}" If $Z$ is a Zariski-closed subvariety of $Y$, then
$$e(Y; x, y) = e(Z; x, y) + e(Y - Z; x, y).$$
So if $Y = \coprod_{i=1}^n Y_i$ is a disjoint union of
finitely many locally closed subsets (i.e. each $Y_i$ is 
the intersection of an open subset and a closed subset), then
$$e(Y; x, y) = \sum_{i=1}^n e(Y_i; x, y).$$
\item"{(2.3)}" If $f: Y \to Z$ is a Zariski-locally trivial bundle 
with fiber $F$, then
$$e(Y; x, y) = e(Z; x, y) \cdot e(F; x, y).$$
\item"{(2.4)}" If $f: Y \to Z$ is a bijective morphism, 
then $e(Y; x, y) = e(Z; x, y)$.
\endroster

As in \cite{Che}, for a complex scheme $Y$, 
we define $e(Y; x, y) = e(Y_{\text{red}}; x, y)$.
\footnote"$\ddag$"{In fact, mixed Hodge structures are defined 
for any complex scheme (see \cite{Del}). Moreover, it is well-known that 
the mixed Hodge structures of $Y$ and $Y_{\text{red}}$
are the same. So a priori we could define the virtual Hodge polynomial
$e(Y; x, y)$ for any complex scheme $Y$. It then follows 
that $e(Y; x, y) = e(Y_{\text{red}}; x, y)$ for any complex scheme $Y$.}
Thus if $Y$ is projective, we have $e(Y; 1, 1) = \chi(Y_{\text{red}})$
by (2.1). It is known that a morphism $f: Y \to Z$ naturally induces 
a morphism $f_{\text{red}}: Y_{\text{red}} \to Z_{\text{red}}$ 
between the reduced schemes $Y_{\text{red}}$ and $Z_{\text{red}}$. 
For instance, $f$ is bijective if and only if $f_{\text{red}}$ is.
Therefore, we see that the above properties 
(2.2)-(2.4) still hold for complex schemes. 

In our applications, the Uhlenbeck compactifications $\frak M^U_{H}(c_1, n)$ 
are complex varieties by J. Li's definition in \cite{LiJ}. 
However, the Gieseker moduli spaces $\frak M^G_{H}(c_1, n)$ are 
complex schemes since they could be nonreduced. 
So the virtual Hodge polynomials $e(\frak M^G_{H}(c_1, n); x, y)$ of 
$\frak M^G_{H}(c_1, n)$ are understood to be 
$e(\frak M^G_{H}(c_1, n)_{\text red}; x, y)$ by our definition.

\section{3. The $S$-duality conjecture in the form formulated 
by Vafa and Witten}

The physics background for the $S$-duality conjecture starts with 
the $N = 4$ supersymmetric Yang-Mills theory which is 
one of the most remarkable known quantum field theories in four dimensions. 
This theory has the largest possible number of supersymmetries 
for a $4$-dimensional theory without gravity. 
A long-standing conjecture asserts that this theory carries 
a symmetry exchanging strong and weak coupling and 
exchanging electric and magnetic fields. 
Montonen and Olive \cite{M-O} proposed a $\Bbb Z_2$ symmetry 
with the above properties and also exchanging the gauge group $G$ 
with its dual group $\hat G$ whose weight lattice is 
the dual of that of $G$ (e.g. $\hat G = SO(3)$ for $G = SU(2)$). 
In fact the $\Bbb Z_2$ symmetry can be extended to a $SL(2, \Bbb Z)$ symmetry, 
known as the $S$-duality, with $A = \pmatrix a&b\cr c&d\cr \endpmatrix$ 
in $SL(2, \Bbb Z)$ acting on the complex parameter 
$\tau = \theta/2\pi + 4\pi i/g^2$ by $A(\tau) = {a\tau + b \over c\tau+d}$ 
where $g$ is the gauge coupling constant and $\theta$ is the theta angle. 
In \cite{V-W}, Vafa and Witten developed a true strong coupling test 
for the $S$-duality conjecture. It is noted that 
the $N=4$ supersymmetric Yang-Mills theory has a twisted version 
that is a topological field theory \cite{Yam}. 
Vafa and Witten showed that its partition function for gauge fields 
in a given topological class is the Euler characteristics of 
instanton moduli spaces when the $4$-dimensional manifold is 
certain algebraic surface such as $\Bbb P^2$ and $K3$. In other words, 
the generating function for the partition functions is 
$$Z_X(\tau, G) = {q^{-s} \over {\# c(G)}} \sum_k \chi(\frak N_k) q^k 
\eqno (3.1)$$
where $q = e^{2\pi i \tau}$, $s$ is some number, 
$X$ is the $4$-dimensional manifold, 
${\# c(G)}$ is the number of elements in the center $c(G)$ of $G$, 
$\chi(\frak N_k)$ is the Euler characteristic of the moduli space $\frak N_k$ 
of $k$-instantons with gauge group $G$. In this setting, 
the $S$-duality conjecture says that $Z_X(\tau, G)$ is modular for 
a finite index subgroup $\Gamma$ of $SL(2, \Bbb Z)$, 
and that for some number $w$, there exists a transformation law
$$Z_X(-1/\tau, G) = \pm \left ( {\tau \over i} \right )^{w/2} 
\cdot Z_X(\tau, \hat G). \eqno (3.2)$$

In the rest of the paper, we fix $G = SU(2)$ and $\hat G = SO(3)$.
Using results about moduli spaces of semistable sheaves 
from algebraic geometry \cite{Got, O'G, Nak, Qi3, Kly, Yos}, 
Vafa and Witten \cite{V-W} discussed the $S$-duality conjecture with 
gauge groups $SU(2)$ and $SO(3)$ for $K3$-surfaces and $\Bbb P^2$. 
It follows from these examples that the numbers $s$ in (3.1) and 
$w$ in (3.2) should be $s = \chi(X)/12$ and $w = -\chi(X)$. 
Thus by (3.1) and (3.2), the generating functions are
$$\align
Z_X(\tau, SU(2)) &= {q^{-\chi(X)/12} \over 2} 
\sum_{k} \chi(\frak N(0, k)) q^{k} \\
Z_X(\tau, SO(3)) &= q^{-\chi(X)/12} 
\sum_{v \in H^2(X, \Bbb Z_2)} \sum_{k} \chi(\frak N(v, k)) q^{k}, \\
\endalign$$
and the transformation law between $Z_X(\tau, SU(2))$ and 
$Z_X(\tau, SO(3))$ is 
$$Z_X(-1/\tau, SU(2)) = \pm 2^{-\chi(X)/2} \left ( {\tau \over i} \right 
)^{- \chi(X)/2} \cdot Z_X(\tau, SO(3))  \eqno (3.3)$$
(note that an additional factor $2^{-\chi(X)/2}$ was added). 
Here $\frak N(0, k)$ is the moduli space of anti-self-dual connections
associated to the $SU(2)$-principal bundle whose second Chern class is $k$,
and $\frak N(v, k)$ with $v \ne 0$ is the moduli space of 
anti-self-dual connections associated to the $SO(3)$-principal bundle 
whose second Stiefel-Whitney class and first Pontryagin class are $v$ and $-4k$
respectively (the instanton number $k$ may not be an integer).
In fact (3.3) has been sharpened. For $v \in H^2(X, \Bbb Z_2)$, let
$$Z_{X, v}(\tau) = q^{-\chi(X)/12} \sum_{k} \chi(\frak N(v, k)) q^{k}. 
\eqno (3.4)$$
So $Z_X(\tau, SU(2)) = {1 \over 2} Z_{X, 0}(\tau)$
and $Z_X(\tau, SO(3)) = \sum_{v \in H^2(X, \Bbb Z_2)} Z_{X, v}(\tau)$.
Then (3.3) is a consequence of the following transformation law:
$$Z_{X, v}(-1/\tau) = \pm 2^{-b_2(X)/2} \left ( {\tau \over i} 
\right )^{- \chi(X)/2} \cdot \sum_{u \in H^2(X, \Bbb Z_2)} 
(-1)^{u \cdot v} Z_{X, u}(\tau) \eqno (3.5)$$
where $(u \cdot v)$ stands for the intersection of $u$ and $v$ modulo $2$. 

Vafa and Witten \cite{V-W} also studied the effect of blowing up 
(or topologically, connect sum with $\overline \Pee^2$) on 
the $S$-duality conjecture. Let $\W X = X \# \overline \Pee^2$
and $E$ be the class of a complex projective line in $\overline \Pee^2$.
Then $H^2(\W X, \Zee_2) \cong H^2(X, \Zee_2) \oplus \Zee_2 E$. 
So a class $\w v$ in $H^2(\W X, \Zee_2)$ can be written as
$\w v = v + a E$ with $v \in H^2(X, \Zee_2)$ and $a = 0$ or $1$. 
Based on (3.5), Vafa and Witten conjectured that
$Z_{\W X, \w v} = Z_a \cdot Z_{X, v}$ where $Z_a$ is a universal function 
independent of $X$ and $v$, i.e. the blowup formula is
$$\sum_{k} \chi(\frak N(\w v, k)) q^{k} = (q^{1 \over 12} \cdot Z_a) \cdot 
\sum_{k} \chi(\frak N(v, k)) q^{k}. \eqno (3.6)$$

We are interested in verifying the blowup formula (3.6) for algebraic surfaces.
For simplicity, we assume that the algebraic surfaces are simply connected.
So let $X$ be a simply connected algebraic surface, 
$\phi: \W X \to X$ be the blowing-up of $X$ at a point $x_0$, 
and $E$ be the exceptional divisor. Fix a divisor $c_1$ on $X$, 
an integer $n \ge [c_1^2/4]$ (the Bogomolov inequality),
and $\w c_1 = \phi^*c_1 - aE$ with $a = 0$ or $1$.
Let $H$ be an ample divisor on $X$ with odd $(H \cdot c_1)$.  
It is well-known \cite{F-M, Bru, Qi1} that for $r \gg 0$, 
all the divisors $H_r = r \cdot \phi^*H - E$ on $\W X$ are ample and  
lie in the same open chamber of type $(\w c_1, n)$.
Thus all the moduli spaces
$\frak M^G_{H_r}(\w c_1, n)$ (resp. $\frak M_{H_r}(\w c_1, n)$) 
with $r \gg 0$ are identical, and shall be denoted by 
$\frak M^G_{H_\infty}(\w c_1, n)$ (resp. $\frak M_{H_\infty}(\w c_1, n)$).
Since $(H_r \cdot \w c_1) = r(H \cdot c_1) - a$ and $(H \cdot c_1)$ is odd,
we can always choose $r \gg 0$ such that $(H_r \cdot \w c_1)$ is also odd.
By a result of Donaldson \cite{Do1}, the Mumford-Takemoto moduli space 
$\frak M_{H}(c_1, n)$ is naturally identified with 
the instanton moduli space $\frak N(v, n-c_1^2/4)$ where $v \equiv c_1 \pmod 2$.
Similarly, $\frak M_{H_\infty}(\w c_1, n)$ is naturally
identified with $\frak N(\w v, n- \w c_1^2/4)$ where 
$\w v \equiv \w c_1 \pmod 2$. Then the blowup formula (3.6) becomes 
$$\sum_{n} \chi(\frak M_{H_\infty}(\w c_1, n)_{\text{red}}) 
  q^{n- {\w c_1^2 \over 4}} = 
(q^{1 \over 12} \cdot Z_a) \cdot 
\sum_{n} \chi(\frak M_{H}(c_1, n)_{\text{red}}) 
  q^{n-{c_1^2 \over 4}}. \eqno (3.7)$$

Since in general $\frak M_{H_\infty}(\w c_1, n)$ and 
$\frak M_{H}(c_1, n)$ are not compact, the blowup formula (3.7)
needs to be modified. First of all, we compactify these moduli spaces. 
There are two compactifications, namely, the moduli spaces 
$\frak M_{H}^G(c_1, n)$ (resp. $\frak M_{H_\infty}^G(\w c_1, n)$) 
of Gieseker semistable sheaves from algebraic geometry 
and the Uhlenbeck compactifications $\frak M^U_{H}(c_1, n)$ 
(resp. $\frak M^U_{H_\infty}(\w c_1, n)$) from gauge theory. 
Unfortunately, these compactification spaces may not be smooth. 
So instead of using the Weil conjecture as in \cite{Yos}, 
we use the virtual Hodge polynomials (Jun Li \cite{LiJ} showed that 
the Uhlenbeck compactifications do carry structures of complex varieties).
Now we finally arrive at the blowup formulae that we shall study 
in the rest of the paper:
$$\sum_{n} e(\frak M^U_{H_\infty}(\w c_1, n); x, y) 
q^{n- {\w c_1^2 \over 4}} = (q^{1 \over 12} \cdot \W Z_a) \cdot 
\sum_{n} e(\frak M^U_{H}(c_1, n); x, y) q^{n- {c_1^2 \over 4}} \eqno (3.8)$$
$$\sum_{n} e(\frak M_{H_\infty}^G(\w c_1, n); x, y) q^{n- {\w c_1^2 \over 4}} 
= (q^{1 \over 12} \cdot \W {\W Z}_a) \cdot 
\sum_{n} e(\frak M_{H}^G(c_1, n); x, y) q^{n- {c_1^2  \over 4}} 
\eqno (3.9)$$
where $\W Z_a = \W Z_a(x, y, q)$ and $\W {\W Z}_a = \W {\W Z}_a(x, y, q)$
are universal functions of $x, y, q, a$.

We shall prove (3.9) in section 5. Formula (3.8) will be proved in 
section 4 by assuming that $\frak M_{H}(c_1, n)$ (respectively,
$\frak M_{H_\infty}(\w c_1, n)$) is dense in $\frak M_{H}^G(c_1, n)$
(respectively, in $\frak M_{H_\infty}^G(\w c_1, n)$) for all $n$.
This condition determines explicitly the structure of 
the Uhlenbeck compactification
$\frak M^U_{H}(c_1, n)$ (respectively, $\frak M_{H_\infty}^U(\w c_1, n)$).
It follows from the formulae in Theorem 4.27 and Theorem 5.27 that 
$$\align
\W Z_a(1, 1, q) &= {\sum_{n \in \Zee} q^{(n+{a \over 2})^2} \over 
  q^{1 \over 12}(1 - q)} \tag 3.10 \\
\W {\W Z}_a(1, 1, q) &= {\sum_{n \in \Zee} q^{(n+{a \over 2})^2} \over 
[q^{1 \over 24} \prod_{n \ge 1} (1 - q^n)]^2}. \tag 3.11 \\
\endalign$$
Following Vafa and Witten \cite{V-W}, 
we rewrite $\W {\W Z}_a(1, 1, q)$ in a different form:
$$\W {\W Z}_a(1, 1, q)={\theta_a(q)\over \eta(q)^2}  \eqno (3.12)$$
where $\theta_a(q)=\sum_{n \in \Zee} q^{(n+{a \over 2})^2}$ and
$\eta(q)=q^{1 \over 24} \prod_{n \ge 1} (1 - q^n)$ is 
the Dedekind $\eta$-function. Our formula (3.12) agrees with 
the conjecture of Vafa-Witten \cite{V-W}.

Finally, we make some remarks about the blowup formula (3.8) and
the Uhlenbeck compactification space. First of all, 
Vafa and Witten \cite{V-W} only used the Gieseker moduli spaces
as the compactification spaces. However, the gauge theoretic compactification
of the moduli space of anti-self-dual instantons is 
the Uhlenbeck compactification space. More importantly, 
the Uhlenbeck compactification works for arbitrary smooth four manifolds 
rather than algebraic surfaces only. So it is natural to
take the Uhlenbeck compactification space into consideration.
Also, as we pointed out earlier, J. Li \cite {LiJ} showed that 
the Uhlenbeck compactification space is a complex algebraic variety. 
So it is an interesting question in its own right 
whether an explicit closed formula (3.8) can be found.

\section{4. Blowup formulae for the Uhlenbeck compactifications}

In this section, we prove the blowup formula (3.8) for 
the Uhlenbeck compactifications. Essentially, the blowup formula (3.8)
is equivalent to a universal relation between the virtual Hodge polynomials
of the Mumford-Takemoto moduli spaces $\frak M_{H_\infty}(\w c_1, n)$
and $\frak M_{H}(c_1, n)$. To prove the universal relation,
we use elementary modifications and stratify $\frak M_{H_\infty}(\w c_1, n)$
into a finite disjoint union of locally closed subsets. 
It turns out that the virtual Hodge polynomials of these subsets
satisfy some universal recursion relations which enable us to 
deduce the universal relation between the virtual Hodge polynomials
of $\frak M_{H_\infty}(\w c_1, n)$ and $\frak M_{H}(c_1, n)$. 
 
First of all, we recall some basic facts about elementary modifications. 
Let $\W V$ be a locally free sheaf on $\W X$ with $c_1(\W V) = \phi^*c_1$.
Then $\W V|_E \cong \Cal O_E(d) \oplus \Cal O_E(-d)$ for some $d \ge 0$.
It is well-known that $d = 0$ if and only if $\W V \cong \phi^*V$
for some locally free sheaf $V$ on $X$. If $d > 0$,
we consider the natural projection 
$\W V \overset \alpha_1 \to \to \Cal O_E(-d) \to 0$.
In fact, since $Hom(\W V, \Cal O_E(-d)) \cong \Cee$,  
the surjection $\alpha_1$ is unique up to scalars.
Let $\W V' = \text{ker}(\alpha_1)$. 
Then $\W V'$ is locally free with $c_1(\W V') = \phi^*c_1 - E$ and 
$c_2(\W V') = c_2(\W V)-d$. Moreover there exists a commutative diagram 
of morphisms
$$\matrix
  &     & 0           &     & 0        &     &              &       \\
  &     & \uparrow    &     & \uparrow &     &              &       \\
0 & \to & \Cal O_E(d) & \to & \W V|_E  & \to & \Cal O_E(-d) & \to 0 \\
  &     & \mapup{\alpha_2}& & \uparrow &     & \Vert        &       \\
0 & \to & \W V'       & \to & \W V  
                 & \overset \alpha_1 \to \to & \Cal O_E(-d) & \to 0 \\
  &     & \uparrow    &     & \uparrow &     &              &       \\
  &     & \W V \otimes \Cal O_{\W X}(-E)    & =   & \W V \otimes \Cal O_{\W X}(-E) &     &              &       \\
  &     & \uparrow    &     & \uparrow &     &              &       \\
  &     & 0           &     & 0        &     &              &       \\
\endmatrix \tag 4.1$$
Conversely, assume that $\W V'$ is a locally free sheaf with 
$c_1(\W V') = \phi^*c_1 - E$ and that there exists a surjection 
$\W V' \overset \alpha_2 \to \to \Cal O_E(d) \to 0$.
Define $\W V$ by putting $\W V \otimes \Cal O_{\W X}(-E) = 
\text{ker}(\alpha_2)$. Then $\W V$ is locally free with 
$c_1(\W V) = \phi^*c_1$ and $c_2(\W V) = c_2(\W V') + d$.
Moreover, this construction also leads to the commutative diagram (4.1).
In fact, the correspondence between the pairs $(\W V, \alpha_1)$ and 
$(\W V', \alpha_2)$ is one-to-one. 

\lemma{4.2} Assume that $(H \cdot c_1)$ is odd. 
Let $\W V$ and $\W V'$ be as in the commutative diagram (4.1). 
Then $\W V$ is $H_r$-stable if and only if $\W V'$ is $H_r$-stable.
\endproclaim
\proof
Since $(H \cdot c_1)$ is odd, $\W V$ is $H_r$-stable
if and only if $(\phi_*\W V)^{**}$ is $H$-stable (see \cite{Bru, Qi1}).
Since $(\phi_*\W V)^{**} \cong (\phi_*(\W V'))^{**}$,
it follows that $\W V$ is $H_r$-stable
if and only if $(\phi_*\W V)^{**}$ is $H$-stable
if and only if $\W V'$ is $H_r$-stable.
\endproof

Next, we study a stratification of $\frak M_{H_\infty}(\w c_1, n)$.
For $d \ge 0$, let $\W {\frak M}_{a, n}(d)$ be the subset of 
$\frak M_{H_\infty}(\w c_1, n)$ parameterizing all locally free sheaves 
$\W V$ with $\W V|_E \cong \Cal O_E(d + a) \oplus \Cal O_E(-d)$.
The subset $\W {\frak M}_{a, n}(d)$ is constructible, i.e.,
a finite disjoint union of locally closed subsets 
(see \cite{Har} for the definition of constructible subsets). 
Then $\W {\frak M}_{0, n}(0) \cong \frak M_{H}(c_1, n)$,
and $\frak M_{H_\infty}(\w c_1, n) = \coprod_{d \ge 0} \W {\frak M}_{a, n}(d)$.
For two nonnegative integers $m_1$ and $m_2$, 
let $U(m_1, m_2)$ be the subset of 
$$\Pee(H^0(\Pee^1, \Cal O_{\Pee^1}(m_1) \oplus \Cal O_{\Pee^1}(m_2))) 
\cong \Pee^{m_1+m_2+1}$$
parameterizing all pairs $(f_1, f_2)$ of homogeneous polynomials 
such that $\deg (f_1) = m_1, \deg (f_2) = m_2$,
and $f_1$ and $f_2$ are coprime. Then $U(m_1, m_2)$ parameterizes 
all surjective maps $\Cal O_{\Pee^1}(-m_1) \oplus \Cal O_{\Pee^1}(-m_2)
\to \Cal O_{\Pee^1} \to 0$. Now the virtual Hodge polynomials of 
$\W {\frak M}_{a, n}(d)$ are given by the following lemma. 

\lemma{4.3} Let $(H \cdot c_1)$ be odd and $\w c_1 = \phi^*c_1 - aE$ 
with $a = 0$ or $1$. Then
\par
\rm{(i)} $\displaystyle{e(\W {\frak M}_{0, n}(d); x, y) 
= \sum_{\ell = 0}^{d-1} e(U(d-\ell -1, d+\ell); x, y)
   e(\W {\frak M}_{1, n-d}(\ell); x, y)}$ for $d > 0$.
\par
\rm{(ii)} $\displaystyle{e(\W {\frak M}_{1, n}(d); x, y) 
= \sum_{\ell = 0}^{d} e(U(d-\ell, d+\ell); x, y)
   e(\W {\frak M}_{0, n-d}(\ell); x, y)}$ for $d \ge 0$.
\endproclaim
\proof 
We shall only prove (i) since similar arguments work for (ii).
For $\W V \in \W {\frak M}_{0, n}(d)$, let $\W V'$ be 
the corresponding elementary modification in (4.1).
Then $\W V' \in \frak M_{H_\infty}(\phi^*c_1 - E, n-d)$.
Let $W_\ell$ be the subset of $\W {\frak M}_{0, n}(d)$ consisting of
all those $\W V$ with $\W V' \in \W {\frak M}_{1, n-d}(\ell)$.
Then $\W {\frak M}_{0, n}(d) = \coprod_\ell W_\ell$. 
Recall that $\W V'|_E \cong \Cal O_E(\ell +1) \oplus \Cal O_E(-\ell)$
for $\W V \in W_\ell$. 
Since $\alpha_2$ in (4.1) is surjective, we must have $\ell \le (d-1)$. So 
$$\W {\frak M}_{0, n}(d) = \coprod_{\ell=0}^{d-1} W_\ell.  \eqno (4.4)$$
\noindent
{\bf Claim:} $e(W_\ell; x, y) = e(U(d-\ell -1, d+\ell); x, y)
\cdot e(\W {\frak M}_{1, n-d}(\ell); x, y)$.
\medskip\noindent
{\it Proof.} By the discussions preceding Lemma 4.2, 
we can prove that $W_{\ell}$ admits a fibration to 
$\W{\frak M}_{1, n-d}(\ell)$ 
with fibers isomorphic to $U(d-\ell-1, d+\ell)$. 
The formula for $e(W_\ell; x, y)$ will follow if 
this fibration is Zariski-locally trivial.
In the following, we shall verify a weaker version 
that there exists a Zariski-locally trivial fibration 
$T_0 \to \W{\frak M}_{1, n-d}(\ell)$ with fibers $U(d-\ell-1, d+\ell)$ 
such that $T_0$ admits a bijective morphism to $W_{\ell}$.
By (2.4), this weaker version is sufficient to prove the Claim.

First, recall that $r \gg 0$ and $(\w c_1 \cdot H_r)$ is odd. 
Thus by the Remark A.7 in \cite{Muk}, there exists a universal bundle 
$\W {\Cal V'}$ over $\W X \times \W {\frak M}_{1, n-d}(\ell)$.
Let $p_i$ be the projection of $\W X\times \W {\frak M}_{1, n-d}(\ell)$ 
to the $i$-th factor, and $\Cal E$ be the relative extension sheaf
$$\Cal Ext^1_{p_2}(p_1^*\Cal O_E(-d),\, \W {\Cal V'}).$$

Next, we show that $\Cal E$ is locally free. Indeed, 
there is a local to global spectral sequence which has $E_2$ term
$E^{m,n}_2= R^mp_{2*}(\Cal Ext^n(p_1^*\Cal O_E(-d), \W{\Cal V}'))$
and converges to $\Cal Ext^{m+n}_{p_2}(p_1^*\Cal O_E(-d), \W{\Cal V}'))$ 
(see \cite{BPS, Fri}). So we have a canonical exact sequence
$$0 \to R^1p_{2*}(\Cal Hom(p_1^*\Cal O_E(-d), \W{\Cal V}')) \to \Cal E \to$$
$$\to p_{2*}(\Cal Ext^1(p_1^*\Cal O_E(-d), \W{\Cal V}')) \to 
R^2p_{2*}(\Cal Hom(p_1^*\Cal O_E(-d), \W{\Cal V}'))$$
Since $p_1^*\Cal O_E(-d)$ is torsion and $\W{\Cal V}'$ is torsion free, 
$\Cal Hom(p_1^*\Cal O_E(-d), \W{\Cal V}') = 0$. Hence 
$$\Cal E \cong p_{2*}(\Cal Ext^1(p_1^*\Cal O_E(-d), \W{\Cal V}'))
\cong p_{2*}(p_1^*\Cal O_E(d-1) \otimes \W{\Cal V}'). 
\eqno (4.5)$$
where we have used the fact that $\Cal Ext^1(\Cal O_E(-d), \Cal O_{\W X})
\cong \Cal O_E(d-1)$. Now 
$$H^0(\Cal O_E(d-1) \otimes \W V') \cong H^0(\Cal O_E(d-\ell-1)\oplus \Cal O_E(d+\ell))$$
for every $\W V' \in \W {\frak M}_{1, n-d}(\ell)$. 
Thus $h^0(\Cal O_E(d-1) \otimes \W V')$
is independent of $\W V' \in \W {\frak M}_{1, n-d}(\ell)$.
So $p_{2*}(p_1^*\Cal O_E(d-1) \otimes \W{\Cal V}')$ is locally free. 
By (4.5), $\Cal E$ is locally free.

Put $T=\Bbb P(\Cal E^*)$. Then $T$ is a Zariski-locally trivial bundle 
over $\W {\frak M}_{1, n-d}(\ell)$.
By (4.5), the fiber of $T$ over $\W V' \in \W {\frak M}_{1, n-d}(\ell)$ 
is canonically isomorphic to 
$$\Pee(Ext^1(\Cal O_E(-d), \W V')) \cong 
\Pee(H^0(\Cal O_E(d-\ell -1)\oplus \Cal O_E(d+\ell))). \eqno (4.6)$$
Moreover, since $Hom(\Cal O_E(-d), \W V')=0$ for 
every $\W V'\in \W {\frak M}_{1, n-d}(\ell)$, 
we obtain a universal extension by the Corollary 4.5 in \cite{Lan}:
$$\exact{q^*_2\Cal L\otimes (\text{Id}_{\W X} \times \gamma)^*\W {\Cal V}'}
{\W {\Cal V}}{q_1^*\Cal O_E(-d)}  \eqno (4.7)$$
where $q_i$ is the projection from $\W X\times T$ to the $i$-th factor, 
$\Cal L$ is the tautological line bundle on $T$, 
and $\gamma$ is the bundle projection
$T\rightarrow \W {\frak M}_{1, n-d}(\ell)$.

Let $T_0$ be the open subset of $T$ consisting of extensions $\xi$:
$$0 \to \W V' \to \W V \to \Cal O_E(-d) \to 0  \eqno (4.8)$$
such that the middle term $\W V$ is locally free. 
By the Proposition 1.30 in \cite{Fri} (see also \cite{Br1, Br2}),
the middle term $\W V$ corresponding to an extension class 
$$\xi \in \Pee(Ext^1(\Cal O_E(-d), \W V')) \subset T$$
is locally free if and if only $\xi \in U(d-\ell -1, d+\ell)$ 
via the canonical identification (4.6). 
Now $p_1^*\Cal O_E(d-1) \otimes \W{\Cal V}'$ is a rank two vector bundle 
over $E \times \W {\frak M}_{1, n-d}(\ell)$, 
and its restriction to $E \times \{ \W V' \}$ with 
$\W V' \in \W {\frak M}_{1, n-d}(\ell)$ is isomorphic to 
$\Cal O_E(d-\ell-1) \oplus \Cal O_E(d+\ell)$. By the Lemma 1 in \cite{Br1},
for each $\W V' \in \W {\frak M}_{1, n-d}(\ell)$,
there exists a Zariski open subset 
$\Cal U \subset \W {\frak M}_{1, n-d}(\ell)$ such that 
$\W V' \in \Cal U$ and 
$$(p_1^*\Cal O_E(d-1) \otimes \W{\Cal V}')|_{E \times \Cal U}
\cong p_1^*(\Cal O_E(d-\ell-1) \oplus \Cal O_E(d+\ell)).$$
Thus via (4.5), $(\gamma|_{T_0})^{-1}(\Cal U)$ is canonically isomorphic to 
$\Cal U \times U(d-\ell -1, d+\ell)$. 
Therefore the projection $\gamma|_{T_0}: T_0 \to \W {\frak M}_{1, n-d}(\ell)$
is Zariski-locally trivial, 
and its fiber over $\W V' \in \W {\frak M}_{1, n-d}(\ell)$ is 
isomorphic to $U(d-\ell -1, d+\ell)$. By (2.3), we have
$$e(T_0; x, y) = e(U(d-\ell -1, d+\ell); x, y)
\cdot e(\W {\frak M}_{1, n-d}(\ell); x, y). \eqno (4.9)$$

We claim that there exists a bijective morphism from $T_0$ to $W_\ell$.
Indeed, by (4.7), (4.8), and Lemma 4.2, 
the bundle $\W {\Cal V}|_{\W X\times T_0}$ induces a morphism
$$\Psi: T_0 \to \frak M_{H_\infty}(\phi^*c_1, n)$$
which sends $\xi \in T_0 \cap \Pee(Ext^1(\Cal O_E(-d), \W V'))$
to the corresponding $\W V$ in (4.8). 
Since $\W V$ is locally free and $c_1(\W V) = \phi^*c_1$,
restricting (4.8) to $E$ yields an exact sequence
$$0 \to \Cal O_E(d) \to \W V|_E \to \Cal O_E(-d) \to 0.$$
Since $d > 0$, the above exact sequence must split, 
i.e. $\W V|_E \cong \Cal O_E(d) \oplus \Cal O_E(-d)$.
Thus by the definition of $W_\ell$, we have $\text{Im}(\Psi) \subset W_\ell$.
To show that $\Psi: T_0 \to W_\ell$ is surjective, 
let $\W V \in W_\ell$ and consider the commutative diagram (4.1).
The surjective map $\alpha_2$ in (4.1) induces a surjective map 
$\Cal O_E(\ell+1) \oplus \Cal O_E(-\ell) \cong \W V'|_E \to \Cal O_E(d) \to 0$, 
i.e., an element $\xi \in U(d-\ell-1, d+\ell)$. 
Regard $\xi \in \Pee(Ext^1(\Cal O_E(-d), \W V'))$ via 
the canonical identification (4.6).
Then the extension determined by $\xi$ is precisely the second row in (4.1).
So $\xi \in T_0$ and $\Psi(\xi) = \W V$. 
Thus $\Psi: T_0 \to W_\ell$ is surjective.
Similarly, using $Hom(\W V, \Cal O_E(-d)) \cong \Bbb C$
for $\W V \in W_\ell$, we can show that $\Psi$ is injective.

Since $\Psi: T_0 \to W_\ell$ is a bijective morphism, 
we obtain by (2.4) and (4.9) that 
$$e(W_\ell; x, y) = e(T_0; x, y) = e(U(d-\ell -1, d+\ell); x, y)
\cdot e(\W {\frak M}_{1, n-d}(\ell); x, y).  \qed$$

Our formula for $e(\W {\frak M}_{0, n}(d); x, y)$ follows from 
(4.4), (2.2), and the Claim.
\endproof

Now using Lemma 4.3, we can prove the universal relation between 
the virtual Hodge polynomials of the moduli spaces
$\frak M_{H_\infty}(\w c_1, n)$ and $\frak M_{H}(c_1, n)$.
In Lemma 4.13 below, we shall partially compute the universal function.
We remark that this universal relation will be used in the proofs of
the blowup formulae for both the Uhlenbeck compactification spaces
and the Gieseker moduli spaces. 

\proposition{4.10} Let $(H \cdot c_1)$ be odd and $\w c_1 = \phi^*c_1 - aE$ 
with $a = 0$ or $1$. Then there exists a universal function $B_a(x, y, q)$ 
such that
$$\sum_{n} e(\frak M_{H_\infty}(\w c_1, n); x, y) q^n
= B_a(x, y, q) \cdot \sum_{n} 
e(\frak M_{H}(c_1, n); x, y) q^n.  \eqno (4.11)$$
\endproclaim
\noindent
{\it Proof.} Note that $\frak M_{H_\infty}(\phi^*c_1, n) = 
\coprod_{d \ge 0} \W {\frak M}_{0, n}(d)$ and 
$\W {\frak M}_{0, n}(0) \cong \frak M_{H}(c_1, n)$.
Using Lemma 4.3 (i) and (ii) iteratively, we conclude that
$$e(\frak M_{H_\infty}(\phi^*c_1, n); x, y) = 
\sum_{k = 0}^{n - [{c_1^2 \over 4}]} B_{0, k}(x, y) \cdot
e(\frak M_{H}(c_1, n -k); x, y) \eqno (4.12)$$
for some universal functions $B_{0, k}(x, y)$. It follows that
$$\sum_{n} e(\frak M_{H_\infty}(\phi^*c_1, n); x, y) q^n
= B_0(x, y, q) \cdot \sum_{n} e(\frak M_{H}(c_1, n); x, y) q^n$$
where we have put $B_0(x, y, q) = \sum_{n \ge 0} B_{0, n}(x, y) q^n$.

Similarly, there exists universal function $B_1(x, y, q)$ such that
$$\sum_{n} e(\frak M_{H_\infty}(\phi^*c_1-E, n); x, y) 
q^n = B_1(x, y, q) \cdot \sum_{n} e(\frak M_{H}(c_1, n); x, y) q^n.  \qed$$

\lemma{4.13} Let $a = 0$ or $1$. Then 
$B_a(1, 1, q) = \sum_{n \in \Zee} q^{n(n+a)}$.
\endproclaim 
\proof
Since a slight change of the proof for $a = 0$ works for $a = 1$,  
we shall only prove the case $a = 0$. 
By the definition of $B_0(x, y, q)$, it suffices to show that
$$B_{0, n}(1, 1) = 
\cases 0, &\text{if $n$ is not a square}\\
       1, &\text{if $n= 0$}             \\
       2, &\text{if $n$ is a square and $n > 0$}. 
\endcases         \eqno (4.14)$$
 From the proof of (4.12), we see that $B_{0, 0}(x, y) = 1$. 
Moreover, when $n > 0$,
$$\align
B_{0, n}(x, y) 
&= \sum_{\Sb 0 \le d_1, 0 \le d_2 \le d_1 -1, \ldots,
   0 \le d_{2s-1} \le d_{2s-2}, 0 \le d_{2s} \le d_{2s-1} -1, \ldots \\
   \sum_{i \ge 1} d_i = n \endSb} \\
&\quad e(U(d_1-d_2 -1, d_1+d_2); x, y) e(U(d_2-d_3, d_2+d_3); x, y) \ldots \\
&\quad e(U(d_{2s-1}-d_{2s} -1, d_{2s-1}+d_{2s}); x, y) \\
&\quad e(U(d_{2s}-d_{2s+1}, d_{2s}+d_{2s+1}); x, y) \ldots. \tag 4.15
\endalign$$
Thus (4.14) holds when $n = 0$. In the following, we assume that $n > 0$.

\noindent
{\bf Claim:} {\it Let $m_1$ and $m_2$ be two integers with 
$0 \le m_1 \le m_2$. Then,
$$e(U(m_1, m_2); 1, 1) =  
\cases 0, &\text{if $m_1 > 0$}\\
       1, &\text{if $m_1 = 0$ and $m_2 > 0$}\\
       2, &\text{if $m_1 = m_2 = 0$}.\\
\endcases  \eqno (4.16)$$
\proof
We use mathematical induction on $m_1$. First of all,
we show that (4.16) is true when $m_1 = 0$. Indeed, 
the subset $U(0, 0)$ of $\Pee(H^0(\Pee^1, 
\Cal O_{\Pee^1}(0) \oplus \Cal O_{\Pee^1}(0))) 
\cong \Pee^{1}$ coincides with $\Pee^1$. 
Since $e(\Pee^d; x, y) = 1 + (xy) + \ldots + (xy)^d$, we have  
$$e(U(0, 0); x, y) = e(\Pee^1; x, y) = 1 + (xy).$$ 
So (4.16) holds for $m_1 = m_2 = 0$. When $m_2 > 0$, the subset $U(0, m_2)$ of 
$$\Pee(H^0(\Pee^1,\Cal O_{\Pee^1}(0) \oplus \Cal O_{\Pee^1}(m_2)))$$
is $\Pee(H^0(\Pee^1,\Cal O_{\Pee^1}(0) \oplus \Cal O_{\Pee^1}(m_2)))-
\Pee(\{ 0 \} \oplus H^0(\Pee^1,\Cal O_{\Pee^1}(m_2))) \cong \Pee^{m_2 + 1} - \Pee^{m_2}$.
Thus,
$$e(U(0, m_2); x, y) = e(\Pee^{m_2 + 1}; x, y) - e(\Pee^{m_2}; x, y)
= (xy)^{m_2+1}.$$ 
Hence (4.16) also holds for $m_1 = 0$ and $m_2 > 0$.

Next let $m_1 > 0$. The possible degree of the greatest common divisor
of a pair 
$$(f_1, f_2) \in \Pee(H^0(\Pee^1, \Cal O_{\Pee^1}(m_1) \oplus 
\Cal O_{\Pee^1}(m_2))) - 
\Pee(\{ 0 \} \oplus H^0(\Pee^1,\Cal O_{\Pee^1}(m_2)))$$
can be $0, \ldots, m_1$. For $d = 0, \ldots, m_1$, let $Y_d$ be the subset of 
$\Pee(H^0(\Pee^1,\Cal O_{\Pee^1}(m_1) \oplus \Cal O_{\Pee^1}(m_2))) 
- \Pee(\{ 0 \} \oplus H^0(\Pee^1,\Cal O_{\Pee^1}(m_2)))$ 
parameterizing all pairs $(f_1, f_2)$ such that the greatest common divisor 
of $f_1$ and $f_2$ has degree $d$. Then we obtain
$$\align
\Pee^{m_1 +m_2+1} - \Pee^{m_2} 
& \cong \Pee(H^0(\Pee^1,\Cal O_{\Pee^1}(m_1) \oplus \Cal O_{\Pee^1}(m_2))) 
- \Pee(\{ 0 \} \oplus H^0(\Pee^1,\Cal O_{\Pee^1}(m_2))) \\
&= \coprod_{d = 0, \ldots, m_1} Y_d. \tag 4.17 
\endalign$$
Let $1 \le d \le m_1$, and $(f_1, f_2) \in Y_d$ with  
$\hbox{gcd}(f_1, f_2) = f$. Then we can write $f_1 = f g_1$ and $f_2 = f g_2$ 
with $f \in \Pee(H^0(\Pee^1,\Cal O_{\Pee^1}(d))) \cong \Pee^d$ and 
$$(g_1, g_2) \in 
\cases U(m_1-d, m_2-d), &\text{if $1 \le d < m_1$}\\
       U(0, m_2-m_1), &\text{if $d = m_1$ and $m_1 < m_2$}\\ 
       U(0, 0) - \{ \text{a point} \}, &\text{if $d = m_1$ and $m_1 = m_2$}.\\
\endcases$$
Thus $Y_d$ is the product of  the space $\Pee^d$ with the space 
 $U(m_1-d, m_2-d)$ when $1 \le d < m_1$ or $d = m_1 < m_2$,
or with the space $U(0, 0) - \{ \text{a point} \} \cong \Pee^1 - \{ \text{a point} \}$ 
when $d = m_1 = m_2$. So for $1 \le d \le m_1$, we have
$$e(Y_d; x, y) = e(\Pee^d; x, y) \cdot 
\cases  
e(U(m_1-d, m_2-d); x, y), &\text{if $1 \le d < m_1$}\\
e(U(0, m_2-m_1); x, y), &\text{if $d = m_1 < m_2$}\\
(xy), &\text{if $d = m_1 = m_2$}\\
\endcases \tag 4.18 $$
Note that $e(\Pee^d; 1, 1) = (d+1)$. By (4.18) and our induction hypothesis,
$$e(Y_d; 1, 1) =  \cases  0, &\text{if $1 \le d < m_1$}\\
(m_1+1), &\text{if $d = m_1$}\\
\endcases \tag 4.19 $$
Since $Y_0 = U(m_1, m_2)$, we conclude from (4.17) and (4.19) that
$$(m_1+ m_2+2) - (m_2+1) = e(U(m_1, m_2); 1, 1) + (m_1+1).  \eqno (4.20)$$
It follows that $e(U(m_1, m_2); 1, 1) = 0$ when $m_1 > 0$.
\endproof

We continue the proof of the formula (4.14). Let $x=y=1$. Then by the Claim,
the typical term in (4.15) is nonzero only if we have 
$$d_1-d_2 -1 =0, d_2-d_3=0, \ldots, d_{2s-1}-d_{2s}-1=0, d_{2s}-d_{2s+1}=0, 
\ldots$$
i.e. $d_{2s} = d_{2s+1} = d_1 - s$. Since $\sum_{i \ge 1} d_i = n$, we obtain
$$n = d_1 + 2(d_1 - 1) + 2(d_1-2) + \ldots + 4 + 2 = d_1^2.  \eqno (4.21)$$ 
It follows that if $n$ is not a square, then $B_{0, n}(1,1) =0$.
If $n$ is a positive square and $n = n_0^2$ with $n_0 > 0$, 
then there is exactly one nonzero term in (4.15) 
$$e(U(0, 2 n_0 -1); 1, 1) e(U(0, 2 n_0 -2); 1,1) \ldots $$
$$\cdot e(U(0, 2 n_0 - (2s-1)); 1, 1) e(U(0, 2 n_0 - 2s); 1, 1) \ldots$$ 
$$\cdot e(U(0, 1); 1, 1) e(U(0, 0); 1,1)$$
given by $d_1 = n_0, d_2 = n_0 -1, \ldots, d_{2s-1} = n_0 -(s-1),
d_{2s} = n_0 - s, \ldots$. Thus by the Claim once again, we conclude that
$B_{0, n}(1,1) = 2$.
\endproof

\par\noindent
{\it Remark 4.22.} (i) Using the expression (4.15), one can compute that
$$\align
&B_{0, 0}(x, y) = 1 \\
&B_{0, 1}(x, y) = (xy)^2(1+xy) \\
&B_{0, 2}(x, y) = (-1) (xy)^2(1+xy) [1 - (xy)^2].\\
\endalign$$
However, it is unclear how to get a closed formula for $B_{0, n}(x, y)$ 
in general.
\par
(ii) Naively $B_{0}(1, 1, q)$ (resp. $B_1(1, 1, q)$) may be obtained from 
the Proposition 0.3 (resp. the Remark 4.5) in \cite{Yos} 
by replacing the $q$ and $t$ there by $1$ and $q$ respectively.
However it is unclear whether an analogue of the Weil conjecture holds 
for arbitrary algebraic varieties and cohomology with compact support. 

\medskip
Let $\frak M^U_H(c_1, n)$ be the Uhlenbeck compactification of
$\frak N(c_1, n-c_1^2/4) \cong \frak M_{H}(c_1, n)$ 
(see \cite{Uhl, Do1, Do2, F-M}), 
and let $\text{Sym}^{n}(X)$ be the $n$-th symmetric product of $X$.
Our next lemma determines explicitly the structure of $\frak M^U_H(c_1, n)$. 

\lemma{4.23} Let $(c_1 \cdot H)$ be odd. Assume that 
the Mumford-Takemoto moduli space $\frak M_{H}(c_1, n)$ 
is dense in the Gieseker moduli space $\frak M_{H}^G(c_1, n)$. Then 
$$\frak M^U_H(c_1, n) = \coprod_{i = [{c_1^2 \over 4}]}^{n} 
\frak M_{H}(c_1, i) \times \Sym^{n-i}(X). \eqno (4.24)$$
\endproclaim
\proof
We follow the argument in the proof of the Proposition 7 in section 5 
of \cite{H-L}. Since $(c_1 \cdot H)$ is odd,
$c_1$ is odd and $H$ does not lie on any wall of type $(c_1, n)$.
Thus, all the relevant anti-self-dual $SO(3)$-connections are irreducible,
and correspond to $H$-stable rank-$2$ bundles
by Donaldson's result \cite{Do1}. It follows from the definition of
the Uhlenbeck compactification that
$$\frak M^U_H(c_1, n) \subset \coprod_{i = [{c_1^2 \over 4}]}^{n} 
\frak M_{H}(c_1, i) \times \Sym^{n-i}(X).$$
Now endow $\frak M_{H}^G(c_1, n)$ with the reduced scheme structure.
By our assumption, $\frak M_{H}(c_1, n)$ is dense in $\frak M_{H}^G(c_1, n)$.
Jun Li's results \cite{LiJ} which are for $c_1=0$ but generalize to
odd $c_1$ say that there exists a morphism
$$\gamma: \frak M_{H}^G(c_1, n) \to \frak M^U_H(c_1, n)$$
such that $\gamma(V) = (V^{**}, Z_V)$,
where $V^{**}$ is the double dual of $V$ and
$$Z_V \quad \overset \text{define} \to = \quad
\sum_{x \in X} h^0(X, (V^{**}/V)_x) \cdot x \in \Sym^{n-c_2(V^{**})}(X).$$

To verify (4.24), it suffices to show that
for $[{c_1^2/4}] \le i \le n$ and for every
$$(V', Z') \in \frak M_{H}(c_1, i) \times \Sym^{n-i}(X),$$
there exists $V \in \frak M_{H}^G(c_1, n)$ such that
$\gamma(V) = (V', Z')$. Choose a $0$-cycle $Z$ on $X$ with
$\sum_{x \in X} \text{length}(Z_x) \cdot x = Z'$. Take any surjective map
$V' \overset \alpha \to \to \Cal O_Z \to 0$. Let $V = \text{ker}(\alpha)$.
Then $V \in \frak M_{H}^G(c_1, n)$ with $V^{**} = V'$ and $Z_V = Z'$.
So $\gamma(V) = (V', Z')$.
\endproof

We remark that our assumption that the Mumford-Takemoto moduli space 
$\frak M_{H}(c_1, n)$ is dense 
in $\frak M_{H}^G(c_1, n)$ is only used in Lemma 4.23 and Theorem 4.27 below.
Moreover, by the Lemma 2.3 in \cite{F-Q}, if $(c_1 \cdot H)$ is odd 
and the anti-canonical divisor $(-K_X)$ is effective, 
then $\frak M_{H}(c_1, n)$ is dense in $\frak M_{H}^G(c_1, n)$.

Let $\frak M^U_{H_\infty}(\w c_1, n)$ be the Uhlenbeck compactification
of $\frak M_{H_\infty}(\w c_1, n)$. Assume that $\frak M_{H_\infty}(c_1, n)$ 
is dense in the Gieseker moduli space $\frak M_{H_\infty}^G(c_1, n)$.
By our convention, we have $\frak M_{H_\infty}(\w c_1, n) \cong 
\frak M_{H_r}(\w c_1, n)$ and $\frak M_{H_\infty}^G(\w c_1, n) \cong 
\frak M_{H_r}^G(\w c_1, n)$ 
where $r \gg 0$ and $(H_r \cdot \w c_1)$ is odd.
Replacing $X, c_1, H$ in (4.24) by $\W X, \w c_1, H_r$ respectively, 
we see that the Uhlenbeck compactification of the moduli space
$\frak M_{H_r}(\w c_1, n)$ is:
$$\coprod_{i = [{\w c_1^2 \over 4}]}^n \frak M_{H_r}(\w c_1, i) 
\times \text{Sym}^{n-i}(\W X)
= \coprod_{i = [{\w c_1^2 \over 4}]}^n \frak M_{H_\infty}(\w c_1, i) 
\times \text{Sym}^{n-i}(\W X).  \eqno (4.25)$$
Therefore the Uhlenbeck compactification $\frak M^U_{H_\infty}(\w c_1, n)$
of $\frak M_{H_\infty}(\w c_1, n)$ is:
$$\frak M^U_{H_\infty}(\w c_1, n) = 
\coprod_{i = [{\w c_1^2 \over 4}]}^n \frak M_{H_\infty}(\w c_1, i) 
\times \text{Sym}^{n-i}(\W X).  \eqno (4.26)$$
Now we are ready to prove the blowup formula (3.8).

\theorem{4.27} Let $(H \cdot c_1)$ be odd and $\w c_1 = \phi^*c_1 - aE$ 
with $a = 0$ or $1$. Assume that $\frak M_{H}(c_1, n)$ 
(respectively, $\frak M_{H_\infty}(\w c_1, n)$) is dense 
in the Gieseker moduli space $\frak M_{H}^G(c_1, n)$ 
(respectively, in $\frak M_{H_\infty}^G(\w c_1, n)$) for all $n$. Then 
$$\sum_{n} e(\frak M^U_{H_\infty}(\w c_1, n); x, y) 
q^{n- {\w c_1^2 \over 4}} = (q^{1 \over 12} \cdot \W Z_a) \cdot 
\sum_{n} e(\frak M^U_{H}(c_1, n); x, y) q^{n- {c_1^2 \over 4}}
\eqno (4.28)$$
where $\W Z_a = \W Z_a(x, y, q)$ is a universal function of $x, y, q, a$ with 
$$\W Z_a(1, 1, q) = {\sum_{n \in \Zee} q^{(n+{a \over 2})^2} \over 
q^{1 \over 12}(1 - q)}. \eqno (4.29)$$
\endproclaim
\noindent
{\it Proof.} Let $h^{s, t}(X)$ stand for the Hodge numbers of $X$. Then
$$\sum_{n} e(\Sym^n(X); x, y) q^n
= {1 \over \prod_{s, t} (1 - x^sy^tq)^{(-1)^{s+t}h^{s, t}(X)}} \eqno (4.30)$$
by the formula (*5) on p.481 of \cite{Che}. By (4.26) and (4.24), we have
$$\sum_{n} e(\frak M^U_{H_\infty}(\w c_1, n); x, y) q^n
= \sum_{n} e(\Sym^n(\W X); x, y) q^n \cdot
\sum_{n} e(\frak M_{H_\infty}(\w c_1, n); x, y) q^n$$
$$\sum_{n} e(\frak M^U_H(c_1, n); x, y) q^n
= \sum_{n} e(\Sym^n(X); x, y) q^n \cdot
\sum_{n} e(\frak M_H(c_1, n); x, y) q^n.$$
Note that $h^{s, t}(\W X) = h^{s, t}(X)$ when $(s, t) \ne (1, 1)$
and $h^{1,1}(\W X) = 1+h^{1,1}(X)$. Therefore by (4.30) and Proposition 4.10, 
we conclude that
$$\align
{\sum_{n} e(\frak M^U_{H_\infty}(\w c_1, n); x, y) q^n \over
   \sum_{n} e(\frak M^U_H(c_1, n); x, y) q^n}
&= {1 \over 1 - xyq} \cdot {\sum_{n} 
    e(\frak M_{H_\infty}(\w c_1, n); x, y) q^n \over 
   \sum_{n} e(\frak M_H(c_1, n); x, y) q^n}   \\
&= {B_a(x, y, q) \over 1 - xyq}.  \tag 4.31 \\
\endalign$$
Since $\w c_1^2 = c_1^2 -a$, the formula (4.28) follows from (4.31) by putting
$$\W Z_a(x, y, q) = {q^{a \over 4} \cdot B_a(x,y, q) \over 
q^{1 \over 12}(1 - xyq)}.  \eqno (4.32)$$
Finally, we obtain (4.29) from (4.32) and Lemma 4.13.
\qed

\section{5. Blowup formulae for the Gieseker moduli spaces}

In this section, we prove the blowup formula (3.9) 
for the Gieseker moduli spaces. By Proposition 4.10, 
it suffices to prove a universal relation between 
the virtual Hodge polynomials of 
the Gieseker moduli space $\frak M^G_{H}(c_1, n)$ and 
the Mumford-Takemoto moduli space $\frak M_{H}(c_1, n)$.
Using standard techniques, we stratify the moduli space $\frak M^G_{H}(c_1, n)$
into a finite disjoint union of locally closed subsets.
We show that these subsets are closely related to $\frak M_{H}(c_1, n-k)$
and the Grothendieck Quot-scheme $\Quot^k_{\Cal O_X^{\oplus 2}}$
where $0 \le k \le n - [c_1^2/4]$. Then we obtain the universal relation 
between the virtual Hodge polynomials of 
$\frak M^G_{H}(c_1, n)$ and $\frak M_{H}(c_1, n)$.

Our first lemma studies the virtual Hodge polynomials of 
Grothendieck Quot-schemes. In its simplest form, 
it says that if $V$ is a locally free rank-$2$ sheaf over $X$, 
then $e(\Quot^n_V; x, y)=e(\Quot^n_{\Cal O_{X}^{\oplus 2}}; x, y)$.
We recall some definitions and notations from \cite{Gro}.
Let $Y$ be a projective scheme over a base noetherian scheme $S$,
and ${\Cal V}$ be a locally free rank-$r$ sheaf over $Y$.
For a nonnegative integer $n$, let $\Quot_{{\Cal V}/Y/S}^n$ be
the (relative) Quot-scheme parameterizing all the surjections
${\Cal V}|_{Y_s} \to Q \to 0$ with $s \in S$ such that the quotients $Q$ are 
torsion sheaves supported at finitely many points and $h^0(Y_s, Q) = n$ 
(for simplicity, we have used $\Quot^n_V$ to stand for 
$\Quot^n_{V/X/\text{Spec}(\Cee)}$). 
Let $\pi: \Quot_{{\Cal V}/Y/S}^n \to S$ be the natural map.
Since $Y \to S$ is projective, so is $\pi$. 
Over $Y \times_S \Quot_{{\Cal V}/Y/S}^n$, there exists a universal quotient:
$$p_1^*\Cal V \to \Cal Q_n \to 0   \eqno (5.1)$$
where $p_1: Y \times_S \Quot_{{\Cal V}/Y/S}^n \to Y$ is the natural projection.

\lemma{5.2} Let ${\Cal V}$ be a locally free rank-$r$ sheaf over $Y$. Then
$$e(\Quot_{{\Cal V}/Y/S}^n; x, y) =
e(\Quot_{\Cal O_Y^{\oplus r}/Y/S}^n; x, y).$$
\endproclaim
\noindent
{\it Proof.} Since ${\Cal V}$ is locally free, we can decompose $Y$ into
a finite disjoint union of locally closed subsets $Y_1, \ldots, Y_m$
such that for each $i$, there exist a Zariski open subset $U_i$ of $Y$
containing $Y_i$ and an isomorphism
${\Cal V}|_{U_i} \cong \Cal O_Y^{\oplus r}|_{U_i}$. 
Let $\Quot_{{\Cal V}/Y/S}^n(Y_i)$ be the subset of
$\Quot_{{\Cal V}/Y/S}^n$ consisting of all the points $q$ such that if
${\Cal V}|_{Y_{\pi(q)}} \to {\Cal Q_n}|_{Y_{\pi(q)}} \to 0$ 
is the surjection parameterized by $q$, then
$$\Supp({\Cal Q_n}|_{Y_{\pi(q)}}) \subset (Y_i)_{\pi(q)}.$$
Since $Y_i$ is locally closed, we conclude that 
the subset $\Quot_{{\Cal V}/Y/S}^n(Y_i)$ is constructible, 
i.e., a finite disjoint union of locally closed subsets. 

First of all, by the universality of Quot-schemes, we claim that
$$\Quot_{{\Cal V}/Y/S}^n(Y_i) \cong
\Quot_{\Cal O_Y^{\oplus r}/Y/S}^n(Y_i)$$
for each $i$. Indeed, restricting (5.1) to
$Y \times_S \Quot_{{\Cal V}/Y/S}^n(Y_i)$ yields a surjection
$$p_1^*\Cal V \overset \beta_i \to \to 
\Cal Q_n|_{Y \times_S \Quot_{{\Cal V}/Y/S}^n(Y_i)} \to 0   \eqno (5.3)$$
over $Y \times_S \Quot_{{\Cal V}/Y/S}^n(Y_i)$
(here and thereafter, by abusing notations, 
we always use $p_1$ to stand for the first projection such as 
$Y \times_S \Quot_{{\Cal V}/Y/S}^n(Y_i) \to Y$).
 From the definition of $\Quot_{{\Cal V}/Y/S}^n(Y_i)$, we see that 
$\Supp(\Cal Q_n|_{Y \times_S \Quot_{{\Cal V}/Y/S}^n(Y_i)})$ 
is contained in $Y_i \times_S \Quot_{{\Cal V}/Y/S}^n(Y_i)$. 
Since $Y_i \subset U_i$, we obtain the inclusions
$$\Supp(\Cal Q_n|_{Y \times_S \Quot_{{\Cal V}/Y/S}^n(Y_i)})
\subset Y_i \times_S \Quot_{{\Cal V}/Y/S}^n(Y_i)    
\subset U_i \times_S \Quot_{{\Cal V}/Y/S}^n(Y_i).  \tag 5.4 $$
Notice that $\text{Hom}(p_1^*\Cal V, 
\Cal Q_n|_{Y \times_S \Quot_{{\Cal V}/Y/S}^n(Y_i)})$ is isomorphic to
$$H^0(Y \times_S \Quot_{{\Cal V}/Y/S}^n(Y_i), (p_1^*\Cal V)^* \otimes
\Cal Q_n|_{Y \times_S \Quot_{{\Cal V}/Y/S}^n(Y_i)}).  \eqno (5.5)$$
By (5.4) and since ${\Cal V}|_{U_i} \cong \Cal O_Y^{\oplus r}|_{U_i}$,
we have an isomorphism 
$$(p_1^*\Cal V)^* \otimes \Cal Q_n|_{Y \times_S \Quot_{{\Cal V}/Y/S}^n(Y_i)} 
\cong (p_1^*\Cal O_Y^{\oplus r})^* \otimes
\Cal Q_n|_{Y \times_S \Quot_{{\Cal V}/Y/S}^n(Y_i)}.  \eqno (5.6)$$
Therefore, we get following isomorphisms
$$\align
&\qquad \text{Hom}(p_1^*\Cal V, 
  \Cal Q_n|_{Y \times_S \Quot_{{\Cal V}/Y/S}^n(Y_i)})\\
&\cong H^0(Y \times_S \Quot_{{\Cal V}/Y/S}^n(Y_i), (p_1^*\Cal V)^* \otimes
  \Cal Q_n|_{Y \times_S \Quot_{{\Cal V}/Y/S}^n(Y_i)})  \\
&\cong H^0(Y \times_S \Quot_{{\Cal V}/Y/S}^n(Y_i), (p_1^*\Cal O_Y^{\oplus r})^* 
  \otimes \Cal Q_n|_{Y \times_S \Quot_{{\Cal V}/Y/S}^n(Y_i)})  \\
&\cong \text{Hom}(p_1^*\Cal O_Y^{\oplus r}, 
  \Cal Q_n|_{Y \times_S \Quot_{{\Cal V}/Y/S}^n(Y_i)}).
\endalign$$
Via these isomorphisms, the map $\beta_i$ in (5.3) induces a map
$$p_1^*\Cal O_Y^{\oplus r} \overset \phi_i \to \to 
\Cal Q_n|_{Y \times_S \Quot_{{\Cal V}/Y/S}^n(Y_i)} $$
over $Y \times_S \Quot_{{\Cal V}/Y/S}^n(Y_i)$. Since $\beta_i$ is surjective, 
we see that the induced map $\phi_i$ is also surjective. 
By the universal property of the Quot-schemes, the surjection 
$$p_1^*\Cal O_Y^{\oplus r} \overset \phi_i \to \to 
\Cal Q_n|_{Y \times_S \Quot_{{\Cal V}/Y/S}^n(Y_i)} \to 0   \eqno (5.7)$$
induces a morphism 
$\Phi_i: \Quot_{{\Cal V}/Y/S}^n(Y_i) \to \Quot_{\Cal O_Y^{\oplus r}/Y/S}^n$.
It is clear that the image of $\Phi_i$ is 
$\Quot_{\Cal O_Y^{\oplus r}/Y/S}^n(Y_i)$ and 
that $\Phi_i$ induces an isomorphism
$$\Quot_{{\Cal V}/Y/S}^n(Y_i) \cong \Quot_{\Cal O_Y^{\oplus r}/Y/S}^n(Y_i).
\eqno (5.8)$$

For every surjection ${\Cal V}|_{Y_s} \to Q \to 0$ parameterized
by a point $q$ in $\Quot_{{\Cal V}/Y/S}^n$ with $\pi(q)=s$,
the quotient $Q$ can be written as
$\oplus_{i=1}^m Q_i$ such that $\Supp(Q_i) \subset (Y_i)_s$.
Moreover, since $Y_1, \ldots, Y_m$ are disjoint,
the surjection ${\Cal V}|_{Y_s} \to Q \to 0$ is equivalent to
the surjections ${\Cal V}|_{Y_s} \to Q_i \to 0$ for $1 \le i \le m$.
It follows from the universal property of
the Quot-schemes that there exists a bijective morphism
$$\coprod_{\sum\limits_i n_i = n} \quad
\prod_{i/S} \Quot_{{\Cal V}/Y/S}^{n_i}(Y_i) \to  
\Quot_{{\Cal V}/Y/S}^n  \eqno (5.9)$$
where we have used the notation
$\prod \limits_{i/S} \Quot_{{\Cal V}/Y/S}^{n_i}(Y_i)$
to stand for the fiber product
of $\Quot_{{\Cal V}/Y/S}^{n_1}(Y_1), \ldots,
\Quot_{{\Cal V}/Y/S}^{n_m}(Y_m)$ over $S$.

More precisely, fix nonnegative integers $n_1, \ldots, n_m$ with
$\sum_i n_i = n$. We shall prove that there exists an injective morphism 
$$\Psi_{n_1, \ldots, n_m}: \prod \limits_{i/S}
\Quot_{{\Cal V}/Y/S}^{n_i}(Y_i) \to \Quot_{{\Cal V}/Y/S}^n  \eqno (5.10)$$
such that $\Quot_{{\Cal V}/Y/S}^n = \coprod\limits_{\sum\limits_i n_i = n}
\text{Im}(\Psi_{n_1, \ldots, n_m})$. Then (5.9) follows immediately. 

To prove (5.10), restricting the universal surjection
$p_1^*\Cal V \to \Cal Q_{n_i} \to 0$ over
$Y \times_S \Quot_{{\Cal V}/Y/S}^{n_i}$ to
$Y \times_S \Quot_{{\Cal V}/Y/S}^{n_i}(Y_i)$, we obtain a surjection
$$p_1^*\Cal V \overset \alpha_i \to \to
\Cal Q_{n_i}|_{Y \times_S \Quot_{{\Cal V}/Y/S}^{n_i}(Y_i)} \to 0$$
over $Y \times_S \Quot_{{\Cal V}/Y/S}^{n_i}(Y_i)$. Consider the map over
$Y \times_S (\prod\limits_{i/S} \Quot_{{\Cal V}/Y/S}^{n_i}(Y_i))$:
$$\alpha \overset \text{def} \to = \W \alpha_1\oplus \ldots\oplus \W \alpha_m:
p_1^*\Cal V \to \bigoplus_i
p_{2, i}^*(\Cal Q_{n_i}|_{Y \times_S \Quot_{{\Cal V}/Y/S}^{n_i}(Y_i)})$$
where $p_{2, i}: Y \times_S (\prod\limits_{i/S}
\Quot_{{\Cal V}/Y/S}^{n_i}(Y_i))
\to Y \times_S \Quot_{{\Cal V}/Y/S}^{n_i}(Y_i)$ is the natural projection,
and $\W \alpha_i: p_1^*\Cal V \to 
p_{2, i}^*(\Cal Q_{n_i}|_{Y \times_S \Quot_{{\Cal V}/Y/S}^{n_i}(Y_i)})$
is the pull-back of $\alpha_i$ via $p_{2, i}$. By (5.4), the support of
$p_{2, i}^*(\Cal Q_{n_i}|_{Y \times_S \Quot_{{\Cal V}/Y/S}^{n_i}(Y_i)})$
is contained in
$$Y_i \times_S (\prod \limits_{i/S} \Quot_{{\Cal V}/Y/S}^{n_i}(Y_i)).$$
Since $Y_1, \ldots, Y_m$ are disjoint and $\alpha_1, \ldots, \alpha_m$
are surjective, $\alpha$ is also surjective.
By the universal property of the Quot-schemes, $\alpha$ induces a morphism
$$\Psi_{n_1, \ldots, n_m}: \prod \limits_{i/S}
\Quot_{{\Cal V}/Y/S}^{n_i}(Y_i) \to \Quot_{{\Cal V}/Y/S}^n.$$
It is clear from the construction that $\Psi_{n_1, \ldots, n_m}$ is injective
and $\Quot_{{\Cal V}/Y/S}^n = \coprod\limits_{\sum\limits_i n_i = n}
\text{Im}(\Psi_{n_1, \ldots, n_m})$. This proves our assertion (5.10). 

Applying (5.10) to $\Cal V$ and $\Cal O_Y^{\oplus r}$, 
and using the identification (5.8) and the properties (2.3) and (2.4) 
of the virtual Hodge polynomials, we get
$$\align  e(\Quot_{{\Cal V}/Y/S}^n; x, y)
&= \sum_{\sum\limits_i n_i = n} e \left (
  \prod \limits_{i/S} \Quot_{{\Cal V}/Y/S}^{n_i}(Y_i); x, y \right )\\
&= \sum_{\sum\limits_i n_i = n} e \left ( \prod\limits_{i/S}
  \Quot_{{\Cal O^{\oplus r}}/Y/S}^{n_i}(Y_i); x, y \right )\\
&= e(\Quot_{{\Cal O^{\oplus r}}/Y/S}^n; x, y).\qed
\endalign$$

Recall that $\W X$ stands for the blowup of the surface $X$ 
at a point $x_0 \in X$ with the exceptional divisor $E$. 
Our next lemma says that there exists a universal relation between
the generating functions 
$\sum_{n} e(\Quot^n_{\Cal O_{X}^{\oplus 2}}; x, y) q^n$ and
$\sum_{n} e(\Quot^n_{\Cal O_{\W X}^{\oplus 2}}; x, y) q^n$.

\lemma{5.11} There exists a universal function $Q(x, y, q)$ such that 
$$\sum_{n} e(\Quot^n_{\Cal O_{\W X}^{\oplus 2}}; x, y) q^n 
= Q(x, y, q) \cdot \sum_{n} 
e(\Quot^n_{\Cal O_{X}^{\oplus 2}}; x, y) q^n. \eqno (5.12)$$
\endproclaim
\proof
Let $\W U = \W X - E$, and let $\W W_k$ be the subset of 
the Grothendieck Quot-scheme $\Quot^n_{\Cal O_{\W X}^{\oplus 2}}$ 
parameterizing all the surjections 
$\Cal O_{\W X}^{\oplus 2} \to \W Q \to 0$ such that 
$\sum_{\w x \in \W U} h^0(\W X, \W Q_{\w x}) = k$. 
Then there exists a bijective morphism
$$\Quot^k_{\Cal O_{\W U}^{\oplus 2}} \times \W T_{n-k} \to \W W_k
\eqno (5.13)$$
where $\W T_{n-k}$ is the subset of $\Quot^{n -k}_{\Cal O_{\W X}^{\oplus 2}}$
parameterizing all the surjections $\Cal O_{\W X}^{\oplus 2} \to \W Q \to 0$
such that the quotients $\W Q$ are supported on $E$.
Since $\Quot^n_{\Cal O_{\W X}^{\oplus 2}} = \coprod_{k = 0}^n \W W_k$,
$$e(\Quot^n_{\Cal O_{\W X}^{\oplus 2}}; x, y) = \sum_{k = 0}^n e(\W W_k; x, y) 
= \sum_{k = 0}^n e(\Quot^k_{\Cal O_{\W U}^{\oplus 2}}; x, y) \cdot 
e(\W T_{n-k}; x, y).$$
So $\displaystyle{\sum_{n} e(\Quot^n_{\Cal O_{\W X}^{\oplus 2}}; x, y) q^n
= \sum_{n} e(\Quot^n_{\Cal O_{\W U}^{\oplus 2}}; x, y) q^n \cdot
\sum_{n} e(\W T_n; x, y) q^n}$. Similarly,
$$\sum_{n} e(\Quot^n_{\Cal O_X^{\oplus 2}}; x, y) q^n
= \sum_{n} e(\Quot^n_{\Cal O_U^{\oplus 2}}; x, y) q^n \cdot
\sum_{n} e(T_{n, x_0}; x, y) q^n$$
where $U = X - \{ x_0 \}$ and $T_{n, x_0}$ is the subset of 
$\Quot^{n}_{\Cal O_X^{\oplus 2}}$ parameterizing all the surjections 
$\Cal O_X^{\oplus 2} \to Q \to 0$ such that $\Supp(Q) = x_0$. 
Since $\W U \cong U$, we obtain
$$\sum_{n} e(\Quot^n_{\Cal O_{\W X}^{\oplus 2}}; x, y) q^n =
{\sum_{n} e(\W T_n; x, y) q^n \over \sum_{n} 
e(T_{n, x_0}; x, y) q^n} \cdot \sum_{n} 
e(\Quot^n_{\Cal O_X^{\oplus 2}}; x, y) q^n. \eqno (5.14)$$

It remains to show that $\W T_n$ and $T_{n, x_0}$
are independent of the surface $X$ and the point $x_0 \in X$.
We shall only prove this for $\W T_n$ since a slight modification 
of the argument also works for $T_{n, x_0}$.
Note that $\W T_n$ is a closed subset of the projective
variety $\Quot^n_{{\Cal O}_{\W X}^{\oplus 2}}$.
By Serre's GAGA principals, it suffices to show that $\W T_n$  
is independent of $X$ and $x_0 \in X$ in analytic category.

Let $B_{x_0}$ be an analytic small open ball containing the point $x_0 \in X$,
and let $\W B_{x_0}$ be the blowup of $B_{x_0}$ at $x_0$.
Since $B_{x_0}$ is independent of $X$ and $x_0$, so is $\W B_{x_0}$.
Now $\phi^{-1}(B_{x_0}) \buildrel\sim\over=\W B_{x_0}$ is 
an open neighborhood of the exceptional divisor $E$. 
Every quotient $\Cal O_{\W X}^{\oplus 2}\rightarrow \W Q\rightarrow 0$
in $\W T_n$ is equivalent to 
the quotient $\Cal O_{\W B_{x_0}}^{\oplus 2}\rightarrow \W Q\rightarrow 0$.
Hence $\W T_n$ is independent of $X$ and $x_0$ in analytic category.
\endproof

Next, we study a stratification of the Gieseker moduli space 
$\frak M^G_H(c_1, n)$. Let $W \in \frak M^G_{H}(c_1, n)$ 
and $W^{**}$ be its double dual. Then $W^{**} \in \frak M_{H}(c_1, n - k)$ 
for some $k$ with $0 \le k \le n - [{c_1^2 \over 4}]$, 
and $W^{**}$ sits in a canonical exact sequence 
$$0 \to W \to W^{**} \to Q \to 0$$
where $Q$ is a torsion sheaf supported at finitely many points
and $h^0(X, Q) = k$. Let
$${\frak M}^k=\{ W\in \frak M^G_{H}(c_1, n)\, |\, 
W^{**} \in \frak M_{H}(c_1, n-k)\}.$$
Then $\frak M^G_H(c_1, n) = \coprod_k {\frak M}^k$ is a decomposition 
of constructible subsets. The following lemma determines 
the virtual Hodge polynomial $e({\frak M}^k; x, y)$ of ${\frak M}^k$. 

\lemma{5.15} Assume that $(H \cdot c_1)$ is odd. Then we have
$$e({\frak M}^k; x, y)= e(\Quot^k_{\Cal O_X^{\oplus 2}}; x, y)
\cdot e(\frak M_{H}(c_1, n-k); x, y).$$
\endproclaim
\noindent
{\it Proof.} Since $(H \cdot c_1)$ is odd, there exists a universal bundle 
$\Cal V$ over $X\times \frak M_{H}(c_1, n-k)$ by the Remark A.7 in \cite{Muk}.
For convenience, we denote $\frak M_{H}(c_1, n-k)$ by $\frak M$ and 
the (relative) Quot-scheme $\Quot^k_{\Cal V/X\times \frak M/\frak M}$ 
by $\Quot$. 

We claim that there exists a bijective morphism $\Psi: \Quot \to {\frak M}^k$.
Indeed, there is a universal surjection
$p_1^*\Cal V \overset \alpha \to \to \Cal Q \to 0$ over 
$(X\times \frak M) \times_{\frak M} \Quot \cong X\times \Quot,$
where $p_1: X\times \Quot \to X \times \frak M$ is the natural projection.
So we have an exact sequence
$$0 \to \text{ker}(\alpha) \to p_1^*\Cal V \overset \alpha \to \to 
\Cal Q \to 0   \eqno (5.16)$$
over $X\times \Quot$. Let $\pi: \Quot \to \frak M$ be the natural projection,
and let $q \in \Quot$. Restrict (5.16) to $X \times q$. 
Since $\Cal Q$ is flat over $\Quot$, we get an exact sequence
$$0 \to \text{ker}(\alpha)|_{X \times q} \to \Cal V|_{X \times \pi(q)} 
\to \Cal Q|_{X \times q} \to 0.$$
Since $\Cal V|_{X \times \pi(q)}$ is $H$-stable and $\Cal Q|_{X \times q}$ 
is a torsion sheaf supported at finitely many points with 
$h^0(\Cal Q|_{X \times q})=k$, we conclude that 
$$\text{ker}(\alpha)|_{X \times q} \in \frak M^k \subset 
\frak M^G_H(c_1, n).$$ 
Hence $\text{ker}(\alpha)$ induces a morphism $\Phi: \Quot \to \frak M^k$.
It is clear that $\Phi$ is bijective. 

By (2.4), $e({\frak M}^k; x, y)= e(\Quot; x, y)$.
Applying Lemma 5.2, we obtain
$$e(\Quot; x, y) =e(\Quot^k_{\Cal O_{X\times \frak M}^{\oplus 2}/
X\times \frak M/\frak M};x, y).$$
By the universal property of Quot-schemes, we have a canonical isomorphism
$$\Quot^k_{\Cal O_{X\times \frak M}^{\oplus 2}/X\times \frak M/\frak M}
\cong \Quot^k_{\Cal O_X^{\oplus 2}/X/\text{Spec}(\Cee)} \times \frak M$$
which is $\Quot^k_{\Cal O_X^{\oplus 2}} \times \frak M$
in our simplified notation. Putting all these together, we get
$$\align
   e({\frak M}^k; x, y)
&= e(\Quot; x, y) = e(\Quot^k_{\Cal O_{X\times \frak M}^{\oplus 2}/
   X\times \frak M/\frak M};x, y) \\
&= e(\Quot^k_{\Cal O_X^{\oplus 2}} \times \frak M; x, y) 
   = e(\Quot^k_{\Cal O_X^{\oplus 2}}; x, y) \cdot e(\frak M; x, y) \\
&= e(\Quot^k_{\Cal O_X^{\oplus 2}}; x, y) \cdot 
   e(\frak M_{H}(c_1, n-k); x, y).  \qed \\
\endalign$$

Now using the universal functions $Q(x, y, q)$ and $B_a(x, y, q)$, 
we prove a universal relation between the virtual Hodge polynomials of 
the Gieseker moduli spaces. In Lemma 5.18 below,
we shall apply this universal relation to determine $Q(1, 1, q)$.  

\proposition{5.17} Let $(H \cdot c_1)$ be odd and $\w c_1 = \phi^*c_1 - aE$ 
with $a = 0$ or $1$. Then 
$$\sum_{n} e(\frak M^G_{H_\infty}(\w c_1, n); x, y) q^n 
= B_a(x, y, q)Q(x, y, q) \cdot \sum_n e(\frak M^G_H(c_1, n); x, y) q^n.$$
\endproclaim 
\noindent
{\it Proof.} Since $\frak M^G_H(c_1, n) = \coprod_k \frak M^k$,
we see from (2.2) and Lemma 5.15 that 
$$\align 
&\qquad  \sum_{n} e(\frak M^G_{H}(c_1, n); x, y) q^n =
  \sum_{n} \sum_k e(\frak M^k; x, y) q^n      \\
&=\sum_{n} \sum_k e(\Quot^k_{\Cal O_X^{\oplus 2}};x,y) 
  \cdot e(\frak M_{H}(c_1, n-k);x,y) \cdot q^n         \\
&=\sum_{n} e(\Quot^n_{\Cal O_{X}^{\oplus 2}}; x, y) q^n \cdot
  \sum_{n} e(\frak M_{H}(c_1, n); x, y) q^n.   \\
\endalign$$
By a similar argument, we also conclude that  
$$\sum_{n} e(\frak M^G_{H_\infty}(\w c_1, n); x, y) q^n = 
\sum_{n} e(\Quot^n_{\Cal O_{\W X}^{\oplus 2}}; x, y) q^n \cdot
\sum_{n} e(\frak M_{H_\infty}(\w c_1, n); x, y) q^n.$$
In view of Proposition 4.10 and Lemma 5.11, we obtain 
$$\sum_{n} e(\frak M^G_{H_\infty}(\w c_1, n); x, y) q^n = B_a(x, y, q)
Q(x, y, q) \cdot \sum_{n} e(\frak M^G_H(c_1, n); x, y) q^n.\qed$$

Next, we shall use the universal relation in Proposition 5.17 and 
the results in \cite{Yos} to determine $Q(1, 1, q)$.
We apply Proposition 5.17 to $X = \Pee^2$. 
However, based on (5.14), it might be possible to determine $Q(1, 1, q)$ 
(resp. $Q(x, y, q)$) by computing $e(\W T_n; 1, 1)$ and $e(T_{n, x_0}; 1, 1)$
(resp. $e(\W T_n; x, y)$ and $e(T_{n, x_0}; x, y)$) directly.

\lemma{5.18} $\displaystyle{Q(1, 1, q) = 
{1 \over \prod_{n \ge 1} (1 - q^n)^2}}$.
\endproclaim 
\par\noindent
{\it Proof.} Let $X = \Pee^2$, $H$ be the divisor represented by a line in $X$, 
$c_1 = -H$, and $a = 0$. Then $c_1 \cdot H = -1$ is odd. 
By the results on p.213 and Theorem 0.4 in \cite{Yos}, 
$$\align 
\sum_n \# \frak M^G_{H}(c_1, n)(\Bbb F_s) q^n
&= \sum_n \# \Quot^n_{\Cal O_X^{\oplus 2}}(\Bbb F_s) q^n \cdot 
\sum_n \# \frak M_{H}(c_1, n)(\Bbb F_s) q^n  \tag 5.19  \\
\sum_n \# \Quot^n_{\Cal O_X^{\oplus 2}}(\Bbb F_s) q^n &= 
\prod_{c \ge 1} \prod_{b=1}^2 Z_s(X, s^{2c-b} q^c)  \tag 5.20 \\
\endalign$$
where $\Bbb F_s$ is a finite field with $s$ elements, 
and $\# Y(\Bbb F_s)$ is the number of rational points 
for an algebraic scheme $Y$ over $\Bbb F_s$, and 
$$Z_s(X, q) \quad {\overset \text{def} \to =} \quad \text{exp} 
\left ( \sum_{r > 0} (\# X(\Bbb F_{s^r})) {q^r \over r} \right )$$
is the zeta function of $X$ over $\Bbb F_s$. Combining (5.19) and (5.20),
we obtain 
$$\sum_n \# \frak M^G_{H}(c_1, n)(\Bbb F_s) q^n
= \prod_{c \ge 1} \prod_{b=1}^2 Z_s(X, s^{2c-b} q^c) \cdot 
\sum_n \# \frak M_{H}(c_1, n)(\Bbb F_s) q^n  \eqno (5.21)$$
By a similar argument, we also conclude that 
$$\sum_n \# \frak M^G_{H_\infty}(\w c_1, n)(\Bbb F_s) q^n
= \prod_{c \ge 1} \prod_{b=1}^2 Z_s(\W X, s^{2c-b} q^c) \cdot 
\sum_n \# \frak M_{H_\infty}(\w c_1, n)(\Bbb F_s) q^n.  \eqno (5.22)$$

Next, by the Proposition 0.3 in \cite{Yos}, we have
$${\sum_n \# \frak M_{H_\infty}(\w c_1, n)(\Bbb F_s) q^n \over 
\sum_n \# \frak M_{H}(c_1, n)(\Bbb F_s) q^n}
= \sum_{n \in \Zee} s^{n(2n-1)} q^{n^2} \cdot \prod_{c \ge 1} 
{1-s^{2c-1}q^c \over 1-s^{2c}q^c}.  \eqno (5.23)$$
Since $Z_s(\W X, q) = {1 \over 1-sq} \cdot Z_s(X, q)$, we obtain
from (5.21), (5.22) and (5.23) that
$${\sum_n \# \frak M^G_{H_\infty}(\w c_1, n)(\Bbb F_s) q^n \over 
\sum_n \# \frak M^G_{H}(c_1, n)(\Bbb F_s) q^n} = 
{\sum_{n \in \Zee} s^{n(2n-1)} q^{n^2} \over \prod_{c \ge 1} (1-s^{2c}q^c)^2}.
\eqno (5.24)$$
It was proved in \cite{Yos} that 
both $\# \frak M^G_{H_\infty}(\w c_1, n)(\Bbb F_s)$ and
$\# \frak M^G_{H}(c_1, n)(\Bbb F_s)$ are polynomials of $s$.
Since all the nonempty moduli spaces $\frak M^G_{H_\infty}(\w c_1, n)$
and $\frak M^G_{H}(c_1, n)$ are smooth, a consequence of the Weil Conjecture 
(see p.197 in \cite{Yos} for more details) says that 
replacing $s$ by $1$ in the right-hand-side of (5.24) yields
$${\sum_n \chi(\frak M^G_{H_\infty}(\w c_1, n)) q^n \over 
\sum_n \chi(\frak M^G_{H}(c_1, n)) q^n} = 
{\sum_{n \in \Zee} q^{n^2} \over \prod_{n \ge 1} (1-q^n)^2}. \eqno (5.25)$$

Since $\frak M^G_{H_\infty}(\w c_1, n)$ is smooth, 
$\chi(\frak M^G_{H_\infty}(\w c_1, n)) = 
e(\frak M^G_{H_\infty}(\w c_1, n); 1, 1)$. Similarly,  
$\chi(\frak M^G_{H}(c_1, n)) = e(\frak M^G_{H}(c_1, n); 1, 1)$. 
Thus by Proposition 5.17 and Lemma 4.13,
$${\sum_n \chi(\frak M^G_{H_\infty}(\w c_1, n)) q^n \over 
\sum_n \chi(\frak M^G_{H}(c_1, n)) q^n} = B_0(1, 1, q)Q(1,1, q)
= \sum_{n \in \Zee} q^{n^2} \cdot Q(1,1, q). \eqno (5.26)$$
It follows immediately from (5.25) and (5.26) that 
$$\displaystyle{Q(1, 1, q) = {1 \over \prod_{n \ge 1} (1 - q^n)^2}}. \qed $$

\theorem{5.27} Let $(H \cdot c_1)$ be odd and $\w c_1 = \phi^*c_1 - aE$ 
with $a = 0$ or $1$. Then
$$\sum_{n} e(\frak M_{H_\infty}^G(\w c_1, n); x, y) q^{n- {\w c_1^2 \over 4}} 
= (q^{1 \over 12} \cdot \W {\W Z}_a) \cdot 
\sum_{n} e(\frak M_{H}^G(c_1, n); x, y) q^{n- {c_1^2  \over 4}}$$
where $\W {\W Z}_a = \W {\W Z}_a(x, y, q)$
is a universal function of $x, y, q, a$ with 
$$\W {\W Z}_a(1, 1, q) = {\sum_{n \in \Zee} q^{(n+{a \over 2})^2} \over 
[q^{1 \over 24} \prod_{n \ge 1} (1 - q^n)]^2}.$$
\endproclaim 
\noindent
{\it Proof.} Follows from Proposition 5.17, Lemma 4.13, and Lemma 5.18. 
Note that
$$\W {\W Z}_a(x, y, q) = {q^{a \over 4} \cdot B_a(x, y, q) \cdot Q(x, y, q)
\over q^{1 \over 12}}.  \qed$$

\enddocument